\documentclass[12pt]{article}

\usepackage{amsmath,amsthm,amssymb}

\newcommand{\1}{\,\mathbf{I}}
\newcommand{\dtp}{\delta_p}
\newcommand{\hh}{\mathcal{H}}
\newcommand{\mtn}{{\rm Mat}_n}
\newcommand{\pp}{\mathcal{P}}

\DeclareMathOperator{\trc}{Tr} 

\begin{document}
\title{A Linear Programming Method for Finding Orthocomplements in Finite Lattices}
\author{George Parfionov\thanks{Friedmann Laboratory For Theoretical Physics,
Department of Mathematics, SPb EF University, Griboyedova 30--32,
191023 St.Petersburg, Russia} and Rom\`an Zapatrin\thanks{Department
of Information Science, The State Russian Museum, In\.zenernaya 4,
191186, St.Petersburg, Russia (corresponding author, e-mail
zapatrin@rusmuseum.ru)}}

\maketitle

\begin{abstract}
A method of embedding partially ordered sets into linear spaces is
presented. The problem of finding all orthocomplementations in a
finite lattice is reduced to a linear programming problem.
\end{abstract}

\emph{MSC 06C15}

\section{INTRODUCTION AND BASIC DEFINITIONS}

We introduce linear algebraic tools for finite lattices. The idea of
the proposed method looks as follows. Given a finite lattice $L$, we
consider its linear hull $H(L)$ as the collection of all real-valued
functions on $L$. Any mapping $f:L\to L$ can be extended to a linear
operator $f:H(L)\to H(L)$. The collection $\hh$ of all linear
operators in the space $H(L)$ is a linear space itself, and we
introdue a convex set of precomplements (see exact definitions
below) as a subset of $\hh$. Then, it turns out that the complements
in the lattice $L$ are in 1--1 correspondence with the solutions of
a linear programming problem in the space $\hh$.

Now we introduce the basic definitions. Given a finite $L$, let
$H=H(L)$ be the set of all real functions on $L$. Fix up a preferred
basis in $H$ labeled by the elements of $L$: for any $p\in L$ its
counterpart is the delta function $\dtp$:
\[
\dtp(q) = \left\lbrace
\begin{array}{cc}
  1 & \mbox{if } p=q \\
  0 & \mbox{otherwise}
\end{array}\right.
\]

The space $H$ possesses the natural structure of commutative algebra
since functions can be pointwise multiplied: for any $f,g\in H$
$(f\cdot g)(q)=f(q)\cdot g(q)$. The unit element $\1$ of $H$ is
$\1=\sum_{p\in L}\dtp$.

\medskip

The partial order $\le$ in $L$ is associated with the zeta operator
\cite{aigner} in $H$, whose matrix is the incidence matrix of the
partial order $L$:
\[
\zeta(p,q) = \left\lbrace
\begin{array}{cc}
  1 & \mbox{if } p\le q \\
  0 & \mbox{otherwise}
\end{array}\right.
\]

The operator $\zeta$ has the following properties: its matrix is
upper triangle (under an appropriate enumeration of the elements of
$L$), and its diagonal entries are equal to 1. Therefore $\zeta$ is
always invertible. Its inverse is denoted by $\mu=\zeta^{-1}$. The
matrix of the operator $\mu$ considered a function of two variables
ranging over $L$ is called the M\"obius function of $L$
\cite{stanley}.

\section{THE POLYTOPE OF PRECOMPLEMENTS}

Two elements $p,q\in L$ are said to be \emph{disjoint} if $p\land q
=0$, \emph{conjoint} if $p\lor q = 1$ (where $0,1$ stand for the
least and the greatest elements of $L$, respectively) and
\emph{complemented} if they are both disjoint and conjoint.A
complement on the lattice $L$ is an idempotent permutation $\alpha$
of the elements of $L$ such that
\begin{itemize}
  \item $p\le q$ implies $\alpha q\,\le\,\alpha_p$
  \item any pair $p$,$\alpha p$ is complemented (it suffices to
  require them to be disjoint)
\end{itemize}

With any permutation $\alpha$ on the lattice $L$ we can associate a
linear operator $\alpha:H\to H$ defined on the basis of delta
functions as follows: $\alpha(\dtp):=\delta_{\alpha p}$. When
$\alpha$ is a complementation on $L$, the matrix of the operator
$\alpha$ is an idempotent orthogonal matrix, that is, satisfying the
following conditions:
\begin{itemize}
  \item $\alpha\ge 0$ --- since the entries of $\alpha$ are 0 or 1
  \item $\alpha\1=\1$ --- since there is only one 1 entry in each
  row
  \item $\alpha^T=\alpha$
\end{itemize}

The next necessary condition for $\alpha$ to be a complement is that
it reverses order. In operator form this is expressed as follows
\[
\alpha\zeta=\zeta^T \alpha
\]
Furthermore, since $\alpha$ is a complement, for any $p,q\in L$ the
conditions $q\le p$ and $q\le \alpha p$ imply $q=0$. This means that
$\zeta\alpha\dtp\cdot\zeta\dtp=\delta_0$, hence
\[
\trc(\zeta^T\zeta\alpha)=\sum_p
\left(\vphantom{I^0_0}\zeta\alpha\dtp,\zeta\dtp\right) =\sum_p
\left(\vphantom{I^0_0}\zeta\alpha\dtp\cdot\zeta\dtp,\1\right)
=\sum_p \left(\vphantom{I^0_0}\delta_0,\1\right) = n
\]
where $p$ ranges over the elements of $L$ and $n$ is the cardinality
of $L$. Now we can introduce the convex subset $\pp$ of the space
$\mtn$ of real $n\times n$ matrices as follows: a matrix $\alpha$ is
in $\pp$ if and only if the following conditions hold:
\begin{equation}\label{eprecomp}
\begin{array}{l}
  \alpha\ge 0 \\
  \alpha=\alpha^T \\
  \trc\left(\vphantom{I^0_0}\zeta^T\zeta\alpha\right)=n \\
  \alpha\1 = \1 \\
  \alpha\zeta=\zeta^T\alpha
\end{array}
\end{equation}
Evidently, all operators $\alpha$ associated with complements are in
$\pp$; however, $\pp$ is a continuous subset of $\mtn$, namely, a
polytope (since all the equations in \eqref{eprecomp} are linear. We
call the polytope $\pp$ defined by \eqref{eprecomp} the
\emph{polytope of pre-complements}.

\section{MAIN RESULTS}

We present two theorems demonstrating the power of the linear
approach to the theory of posets.

\medskip

\emph{Theorem 1.} The orthocomplements, and only they, are the
integer vertices of the polytope $\pp$.

\emph{The idea of the proof.} Otherwise the condition $\alpha\1=\1$
will be broken.

\medskip

\emph{Theorem 2.} The orthocomplements of the lattice $L$ are the
optimal solutions of the following linear programming problem:
\[
\begin{array}{l}
  \trc\left(\vphantom{I^0_0}\zeta\zeta^T\alpha\right) \rightarrow \min \\
  \left\lbrace\begin{array}{ccc}
                \alpha\1 & = & \1 \\
                \alpha^T & = & \alpha \\
                \alpha\zeta & = & \zeta^T\alpha
              \end{array}
              \right.
\end{array}
\]

\emph{The idea of the proof.} The condition $\alpha\1=\1$ implies
$\trc\left(\vphantom{I^0_0}\zeta\zeta^T\alpha\right)\ge n$, then use
the previous theorem.

\section{CONCLUDING REMARKS}

The proposed techniques of linear embedding of finite lattices can
be extended to posets. All the results remain valid, but we face a
strange effect: in the linear hull $H$ of \emph{any} poset $L$ the
operations $\land$ and $\lor$ are always well-defined, but the
results of these operations may bring us beyond the poset $L$ in
question and yield a weighted sum of the elements of the poset. This
issue needs further investigation.

\end{document}